\documentclass{amsart}

\def\Q{{\mathbf Q}}
\def\Z{{\mathbf Z}}

\def\Gal{\mathrm{Gal}}
\def\End{\mathrm{End}}
\def\Aut{\mathrm{Aut}}
\def\Hom{\mathrm{Hom}}
\def\I{{I}}
\def\J{{J}}
\def\F{{\mathbf F}}

\def\GL{\mathrm{GL}}

\def\Sp{\mathrm{Sp}}

\def\M{\mathrm{M}}
\def\dim{\mathrm{dim}}
\def\O{{\mathcal O}}
\def\P{{P}}

\def\X{{\mathcal X}}
\def\Y{{\mathcal Y}}
\def\T{{\mathcal T}}
\def\B{{B}}
\def\invlim{{\displaystyle{\lim_{\leftarrow}}}}

\newtheorem{thm}{Theorem}[section]
\newtheorem{lem}[thm]{Lemma}
\newtheorem{cor}[thm]{Corollary}
\newtheorem{prop}[thm]{Proposition}
\theoremstyle{definition}
\newtheorem{defn}[thm]{Definition}
\newtheorem{asmpt}[thm]{Assumption}
\newtheorem{ex}[thm]{Example}
\newtheorem{rem}[thm]{Remark}

\title[\'Etale cohomology]
{\'Etale cohomology and reduction of abelian varieties}

\author[A.\ Silverberg]{A.\ Silverberg}
\address{Mathematics Department, Ohio State University, 
Columbus, Ohio, 
USA}
\email{silver\char`\@math.ohio-state.edu}

\author[Yu. G. Zarhin]{Yu. G. Zarhin}
\address{Mathematics Department, Pennsylvania State University, 
University Park, PA, 
USA\newline\indent
Institute for Mathematical Problems in Biology, Russian
Academy of Sciences, Push\-chi\-no, Moscow Region, 
Russia}
\email{zarhin\char`\@math.psu.edu}

\begin{document}

\begin{abstract}
In this paper we study the \'etale cohomology groups
associated to abelian varieties. We obtain necessary
and sufficient conditions for an abelian variety to have
semistable reduction (or purely additive reduction
which becomes semistable over a quadratic extension) 
in terms of the action of the
absolute inertia group on the \'etale 
cohomology groups with finite coefficients.
\end{abstract}

\maketitle

\section{Introduction}
Suppose $X$ is a smooth projective variety over a field $F$,
$v$ is a discrete valuation on $F$, and $\ell$ is
a prime number not equal to the residue characteristic of $v$.
Let $F^s$ denote a separable closure of $F$, let ${\bar v}$
be an extension of $v$ to $F^s$, let $\I$ denote the
inertia subgroup at ${\bar v}$ of $\Gal(F^s/F)$, and let
${\bar X} = X\times_F F^s$. For every
positive integer $k$, the group $\I$ acts naturally
on the $k$-th $\ell$-adic \'etale cohomology group 
$H^k_{\text{\'et}}({\bar X},\Q_\ell)$.
Grothendieck proved the Monodromy Theorem (see the Appendix to
\cite{SerreTate}, and 1.2 and 1.3 of
\cite{SGADeligne}), which says that $\I$ acts 
on $H^k_{\text{\'et}}({\bar X}, \Q_\ell)$ 
via quasi-unipotent operators, i.e.,
for every $\sigma\in\I$ we have 
$(\sigma^m-1)^{r} H^k_{\text{\'et}}({\bar X},\Q_\ell)
=0$
for some positive integers $m$ and $r$.
It is known (see 3.7 of \cite{SGADeligne}, and 
3.5 and 3.6 of \cite{SGA7}) 
that if $k=1$, then one may take $r=2$. 
It easily follows (see 
Theorem \ref{cohcor}i,ii below) that if 
$X$ is an abelian variety, then one may take $r=k+1$.
It is shown in \cite{SGA7I} (see 3.4 and 3.8 of \cite{SGADeligne},
and p.~VI of \cite{SGAintro})
that one may take $r=k+1$ whenever one knows the Purity
Conjecture (3.1 of \cite{SGADeligne}) and resolution of singularities.

 From now on, suppose $X$ is a $d$-dimensional abelian variety.
The N\'eron model $\X$ of $X$ at $v$ is a smooth 
separated model of $X$ over the valuation ring $R$
such that for every smooth scheme $\Y$ over $R$ and morphism 
$\varphi : \Y \otimes_R F \to X$ over $F$ there is a unique morphism
$\Y \to \X$ over $R$ which extends $\varphi$. 
The generic fiber of $\X$
can be canonically identified with $X$, and $\X$ is a commutative group
scheme over $R$ whose group structure extends that of $X$.
Let $X_v^0$ denote the identity component of 
the special fiber of $\X$ at $v$. 
Over an algebraic closure of the residue field, there is
an exact sequence of algebraic groups 
$$0 \to U \times \T \to X_{v}^{0} \to \B \to 0,$$
where $\B$ is an abelian variety, $\T$ is the maximal
algebraic torus in $X_{v}^{0}$, 
and $U$ is a unipotent group.
By definition, $X$ is {\it semistable} at $v$ 
if and only if $U=0$.
As $\I$-modules,
$H^1_{\text{\'et}}({\bar X}, \Z_\ell)$
and the $\ell$-adic Tate module $T_{\ell}(X)$ are isomorphic.
Grothendieck's Galois Criterion for Semistability  
says that $X$ is semistable at $v$ if and only if
every $\sigma\in\I$ acts on  $T_{\ell}(X)$ 
as a unipotent operator of echelon $\le 2$, i.e.,
if and only if
$(\sigma-1)^2H^1_{\text{\'et}}({\bar X}, \Z_\ell)=0$
for every $\sigma\in\I$.

Suppose $n$ is a positive integer relatively prime to the
residue characteristic of $v$. 
Then $H^1_{\mathrm{\acute{e}t}}(\bar{X},\Z/n\Z)$,
viewed as an $\I$-module, is isomorphic
to the $n$-torsion $X_n$ on $X$.
Raynaud's criterion says that if 
$\I$ acts trivially on 
$X_n$, and $n \ge 3$,
then $X$ is semistable at $v$.
The authors 
(see \cite{degree} and \cite{banff})  proved that
if $n \ge 5$, then
$X$ is semistable at $v$ if and only if
$(\sigma-1)^{2} X_n=0$
for every $\sigma \in \I$.
In other words, necessary and sufficient conditions for 
semistability can be read off not only
from the $\ell$-adic representation, as shown by
Grothendieck, but also from the mod $n$ representation
(for $n \ge 5$).
The aim of this paper 
(see Theorem \ref{highercohcor})
is to generalize this result to the 
case of the higher \'etale cohomology groups $H^k$. 
Assume that $0<k<2d$, that $k < r\in\Z$, and that 
$n$ does not belong
to a certain finite set  $N(r)$ of prime powers, 
defined explicitly in terms of $r$ in \S\ref{notation}.
(For example, $N(2)=\{1,2,3,4\}$.)  We show that if $k$ is odd, 
then $X$ is semistable at $v$ if and only if
$(\sigma-1)^{r} H^{k}_{\text{\'et}}({\bar X}, \Z/n\Z)
=0$
for every $\sigma \in \I$.
If $k$ is even, we show  
(under an additional assumption; see \S\ref{rmksect}) that
$(\sigma-1)^{r} H^{k}_{\text{\'et}}({\bar X}, \Z/n\Z)
=0$
for every $\sigma \in \I$
if and only if either $X$ is semistable 
at $v$ or
$X$ has purely additive reduction at $v$ but is 
semistable over a (ramified) quadratic extension of $F$.

We treated the case $k=1$ (and $r=2$)
in \cite{degree} and \cite{banff}. 
By Poincar\'e duality one can therefore
treat the case $k=2d-1$. One can check that in the cases
$k=1$ or $2d-1$,  in the criteria above 
one cannot replace $N(r)$ by a smaller set. 
However, when $2 \le k \le 2d-2$ 
one may replace $N(r)$ by an explicitly defined subset $N'(r)$
for which the result is sharp (see \S\ref{midk}). 

In \S\ref{notation} we
introduce basic definitions and notation. 
Section \ref{linalgsect} deals with multilinear algebra
in characteristics $0$ and $\ell$ and over $\Z_\ell$.
We use the Jordan decompositions of exterior powers
of linear operators to obtain a Minkowski-Serre type result.
Section \ref{lemsect} contains
some abelian variety results that will be used later. 
In \S\S\ref{highcohsect}-\ref{midk} 
we state and prove our main results.
We give necessary and sufficient
conditions for semistability, and also necessary and
sufficient conditions for an abelian variety 
to either be semistable or have
purely additive reduction which becomes semistable over a
quadratic extension. 
In \S\ref{midk} we shrink the exceptional set 
in the criteria when
$2 \le k \le 2d-2$. 
We prove that this exceptional set is minimal.

We hope that our results and/or methods will be useful in the study of
semistability for the more general class of motives
\cite{Motives}.

Silverberg would like to thank IHES, the Bunting
Institute, and the Mathematics Institute of the University of
Erlangen-N\"urnberg for their hospitality, 
and NSA, NSF, the Science Scholars
Fellowship Program at the Bunting Institute,
and the Alexander von Humboldt-Stiftung for financial support.
Zarhin would like to thank NSF and Universit\'e Paris-Sud
for financial support and
Centre Emil Borel (Institut Henri Poincar\'e)
and University of Glasgow for their hospitality.
He also would like to thank Ed Formanek for helpful discussions.

\section{Notation and definitions}
\label{notation}

If $F$ is a field, let $F^s$ denote a separable closure. 
Throughout this paper, 
$X$ is a $d$-dimensional abelian variety defined over $F$, and 
$v$ is a discrete valuation on $F$ of residue characteristic $p \ge 0$.
Let ${\bar v}$ denote an extension of $v$ to $F^s$, and 
let $\I$ denote the inertia subgroup at ${\bar v}$ of $\Gal(F^s/F)$.
If $\ell$ is a prime not equal to the characteristic of $F$,
let
$$\rho_{\ell,X} : \Gal(F^s/F) \to \Aut(T_\ell(X)) \cong \GL_{2d}(\Z_\ell)$$
denote the $\ell$-adic representation on the Tate module $T_\ell(X)
= \invlim X_{\ell^n}$. 
Let $V_{\ell}(X)=T_\ell(X)\otimes_{\Z_{\ell}}\Q_{\ell}$,
let ${\bar X} = X \times_F F^s$, and
let  
$$H^k_\ell = H^k_{\text{\'et}}({\bar X},\Z_\ell)
=\Hom_{\Z_{\ell}}(\wedge^k(T_\ell(X)),\Z_\ell).$$
Then $H^k_\ell$ is a free $\Z_\ell$-module of rank $\binom{2d}{k}=:b$, and
$H^k_\ell \otimes_{\Z_\ell}\Q_\ell = \wedge^k_{\Q_{\ell}}(V^\ast)$,
where $V^\ast=\Hom_{\Q_{\ell}}(V,\Q_{\ell})$ is the dual vector space. 
Let 
$$\rho_{\ell,k} : \I \to \Aut(H^k_\ell) \cong \GL_b(\Z_\ell)$$
be the representation giving the action of $\I$ on $H^k_\ell$.

\begin{defn}
If $r$ is a positive integer, define a finite set 
of prime powers $N(r)$ by 
$$N(r)=
 \{\text{prime powers $\ell^m \mid 0 \le m(\ell - 1) \le r $}\}.$$
\end{defn}
For example, $N(1) = \{1, 2\}$, $N(2) = \{1, 2, 3, 4\}$,
$N(3) = \{1, 2, 3, 4, 8\}$,    
$N(4) = \{1, 2, 3, 4, 5, 8, 9, 16\}$.


\section{Some linear algebra}
\label{linalgsect}

\begin{lem}
\label{level2lem}
Suppose $g$ is a linear operator on a finite-dimensional
vector space over a field of characteristic zero. 
Suppose $(g^m-1)^2=0$ for some positive
integer $m$. If $g$ is unipotent, then
$(g - 1)^2 = 0$. If $-g$ is unipotent, then
$(g + 1)^2 = 0$.
\end{lem}

\begin{proof}
If $g$ is unipotent then
all the eigenvalues of $g^{m-1}+ \cdots +g+1$ are $m \ne 0$. 
Therefore $g^{m-1}+ \cdots +g+1$ is an invertible operator.
Since $$0=(g^m-1)^2=(g-1)^2(g^{m-1}+ \cdots +g+1)^2$$
we have $(g-1)^2=0$.

If $-g$ is unipotent then all the eigenvalues of 
$g^m$ are $(-1)^m$. Since $g^m$ is unipotent, $m=2r$ is even. 
We have
$$0=(g^m-1)^2=(g^{2r}-1)^2=
(g+1)^2(g-1)^2(g^{2(r-1)}+ \cdots +g^2+1)^2.$$
All the eigenvalues of $g-1$ are $-2$ and all the 
eigenvalues of $g^{2(r-1)}+ \cdots +g^2+1$ are $r \ne 0$. 
Therefore $g-1$ and $g^{2(r-1)}+ \cdots +g^2+1$ are invertible,  
and $(g+1)^2=0$.
\end{proof}

\begin{thm}
\label{localglobal}
Suppose 
$m$, $r$ and $g$ are positive integers, 
$\ell$ is a prime number, and  
$A \in \M_g(\Z_{\ell})$ satisfies 
$(A - 1)^r \in \ell^m\M_g(\Z_{\ell})$.
Suppose $\lambda$ is an eigenvalue of $A$ which is
a root of unity.
If $r= m(\ell-1)$ then $\lambda^{\ell}=1$.
If $r< m(\ell-1)$ then $\lambda=1$.
\end{thm}

\begin{proof}
This is a special case of Theorem 6.7 of \cite{serrelem}.
\end{proof}

\begin{thm}
\label{localglobal2}
Suppose ${\bar \Z}$ is the ring of algebraic integers in $\bar{\Q}$,
$n$ and $r$ are positive integers, $\lambda$ is
a root of unity in ${\bar \Z}$, and 
$(\lambda-1)^r \in n{\bar \Z}.$
Then  
$\lambda = 1$ if $n \notin N(r)$.
\end{thm}

\begin{proof}
This is a special case of Corollary 3.3 of \cite{serrelem}.
\end{proof}

\begin{lem}[Lemma 4.3 of \cite{degree}]
\label{groupring} 
Suppose $\O$ is an integral domain of characteristic zero, and
$\ell$ is a prime number. Suppose $r$, $s$, and $m$ are positive integers
such that $r \ge m\ell^{s-1}(\ell-1)$. Suppose $\alpha \in \O$ and
$\alpha^{\ell^s} = 1$. Then 
$(\alpha-1)^{r} \in \ell^m\Z[\alpha]$.
\end{lem}

Suppose $V$ is a finite-dimensional vector space over $\Q_\ell$.
Suppose $k$ is an integer and $0 < k < \dim(V)$.
Let 
$$f_{k} : \GL(V) \to \GL(\wedge^k(V))$$
denote the natural representation defined by $f_k(g)=\wedge^k(g)$.
For $g \in \GL(V)$, let 
$g=s_{g}u_{g} = u_{g}s_{g}$ 
be the Jordan decomposition,
where $u_{g}$, $s_{g} \in\GL(V)$, 
$u_{g}$ is unipotent, and
$s_{g}$ is semisimple.
Then $
f_k(g) =
f_k(s_{g})f_k(u_{g}) = f_k(u_{g})f_k(s_{g})$,
$f_k(u_{g})$ is unipotent, and $f_k(s_{g})$ is semisimple. 

\begin{rem}
\label{astrem}
In the notation of \S\ref{notation}, with $^\ast$ denoting 
the dual map, for every $\sigma\in\I$ we have 
$$\rho_{\ell,k}(\sigma)=(f_{k}(\rho_{\ell,X}(\sigma^{-1})))^\ast.$$
\end{rem}

\begin{rem}
\label{unipzl0}
The kernel  of $f_{k}$ is
 $\{\gamma \in \Z_\ell \mid \gamma^k=1\}$.
In particular, if $(k,\ell-1)=1$ then  $f_{k}$ is injective.
\end{rem}

\begin{lem}
\label{unipzl}
The element $\wedge^k(g)$ is unipotent if and only if
there exists a $k$-th root of unity $\gamma \in \Z_\ell$ such that 
$\gamma g$ is unipotent.
\end{lem}

\begin{proof}
If $\wedge^k(g)$ is unipotent, then  $\wedge^k(s_g)=1$.
By Remark \ref{unipzl0}, $s_g$ is a
$k$-th root of unity in $\Z_\ell$. Let $\gamma=s_g^{-1}$. 
\end{proof}

Now assume that $V$ is even-dimensional, choose a non-degenerate  
alternating bilinear form on $V$, and let $\Sp(V) \subset \GL(V)$ be 
the corresponding symplectic group. Write
$$\rho_k:\Sp(V) \to  \GL(\wedge^k(V))$$
for the restriction of $f_k$ to $\Sp(V)$. It is well-known that
$u_g, s_g \in \Sp(V)$ for every $g \in \Sp(V)$.

\begin{rem}
\label{lini}
Remark \ref{unipzl0} easily implies that the kernel of $\rho_{k}$
is $\{1\}$ if $k$ is odd and is $\{1, -1\}$ if $k$ is even.  
\end{rem}

\begin{lem}
\label{lin}
Suppose that $g \in \Sp(V)$.
\begin{enumerate}
\item[{(i)}] If $(g-1)^{2}=0$, then
$(\rho_{k}(g) - 1)^{k+1} = 0$.
\item[{(ii)}] If $k$ is even and $(g+1)^{2}=0$, then
$(\rho_{k}(g) - 1)^{k+1} = 0$.
\item[{(iii)}] 
If $\gamma g$ is unipotent for some $\gamma \in
\Q_\ell$, then $\gamma\in\{\pm 1\}$.
\item[{(iv)}]
Suppose $\wedge^k(g)$ is unipotent. 
If $k$ is odd then $g$ is unipotent.
If $k$ is even then either $g$ or $-g$ is unipotent.
\end{enumerate} 
\end{lem}

\begin{proof}
Note that 
$$(\rho_{k}(g)-1)(v_{1}\wedge v_{2}\ldots \wedge v_{k})
= gv_{1}\wedge\ldots\wedge gv_{k} - 
v_{1}\wedge v_{2}\ldots \wedge v_{k}.$$
Part (i) follows by substituting $g=1+\eta$, and
(ii) follows from (i) applied to $-g$.
Suppose now that $g$ is not unipotent. 
Then $g$ has an eigenvalue $\lambda \ne 1$. Since $g \in \Sp(V)$,
$\lambda^{-1}$ is also an eigenvalue of $g$. 
Suppose also that $\gamma g$ is unipotent for some $\gamma \in
\Q_\ell$.
Then $\gamma\lambda$ and $\gamma/\lambda$ are eigenvalues of the
unipotent element $\gamma g$ and
thus are equal to $1$. Therefore 
$\gamma^2=1$,
i.e., $\gamma\in\{\pm 1\}$, so either $g$ or $-g$ is unipotent. 
If $\wedge^k(g)$ is unipotent, then 
$\gamma g$ is unipotent for some 
$\gamma\in\Q_\ell$, by Lemma \ref{unipzl}.
If $k$ is odd then $\wedge^k(-g)=-\wedge^k(g)$ is not unipotent, so
$-g$ is not unipotent and thus $g$ is unipotent.
\end{proof}

Theorem \ref{unipex} below will be used in \S\ref{midk}.
To prove it, we first prove a lemma and a theorem.

\begin{lem}
\label{newlinlem}
Suppose $\mathcal{W}$ is an $(\ell-1)$-dimensional
vector space over a field of prime characteristic
$\ell \ge 5$, and 
$A$ is a unipotent linear operator on $\mathcal{W}$ 
whose Jordan form
consists of one Jordan block of size $\ell-1$.
 If $r\in \Z$ and
$2 \le r \le \ell-3$, then
$(\wedge^r(A)-1)^{\ell-1} \ne 0$.
\end{lem}

\begin{proof} 
It follows immediately
from Corollary III.2.7(a) on p.~43 and
Proposition III.2.10(b) on p.~45 of \cite{AF}
that all but one of the Jordan blocks of $\wedge^r(A)$
have size $\ell$, and the size of the remaining block
is $1$ or $\ell-1$ (since the size is less than
$\ell$ and is congruent mod $\ell$ to
$\dim(\wedge^r(\mathcal{W})) = {\binom{\ell - 1}{r}}
\equiv (-1)^r$).
Since $\ell\ge 5$ and $2 \le r \le \ell-3$, we have
$\dim(\wedge^r(\mathcal{W}))\ge\ell$. Therefore,
$\wedge^r(A)$
must have a Jordan block of size $\ell$, so
$(\wedge^r(A)-1)^{\ell-1} \ne 0$.
\end{proof}

\begin{thm}
\label{unip}
Suppose ${\mathcal V}$
is a finite-dimensional vector space over a field of 
characteristic $\ell \ge 5$, and
$A$ is a unipotent linear operator on ${\mathcal V}$. 
Suppose $k \in \Z$, 
 $2 \le k \le \dim({\mathcal V})-2$, and
$(A-1)^{\ell-2} \ne 0$.
Then $(\wedge^k(A)-1)^{\ell-1} \ne 0$.
\end{thm}

\begin{proof} 
Let $\mathcal{W}$ be an $A$-invariant 
$(\ell-1)$-dimensional subspace of 
${\mathcal V}$ such that the Jordan form of the
restriction of $A$ to $\mathcal{W}$ 
is a Jordan block of size $\ell-1$.
Then $\wedge^j(\mathcal{W})$ 
is a $\wedge^j(A)$-invariant subspace of
$\wedge^j({\mathcal V})$, for all $j$. 
By Lemma \ref{newlinlem}, 
$$(\wedge^{j}(A)-1)^{\ell-1}(\wedge^{j}(\mathcal{W})) \ne 0$$
if $2 \le j \le \ell-3$. We are therefore done if 
$k \le \ell-3$.
Suppose $k>\ell-3$. Then 
$\wedge^{k-(\ell-3)}(\mathcal{V}/\mathcal{W}) \ne 0$.
Since $A$ is unipotent, $\wedge^{k-(\ell-3)}(A)$ is also unipotent, 
and there exists a non-zero $\wedge^{k-(\ell-3)}(A)$-invariant
element $u\in\wedge^{k-(\ell-3)}({\mathcal V}/\mathcal{W})$.
Then
$$0 \ne (\wedge^{\ell-3}(A)-1)^{\ell-1}(\wedge^{\ell-3}(\mathcal{W})) \cong
(\wedge^{\ell-3}(A)-1)^{\ell-1}(\wedge^{\ell-3}(\mathcal{W}))\otimes u$$ 
$$ = 
(\wedge^{k}(A)-1)^{\ell-1}(\wedge^{\ell-3}(\mathcal{W})\otimes u)
$$
$$\subset  
(\wedge^{k}(A)-1)^{\ell-1}(\wedge^{\ell-3}(\mathcal{W})\otimes 
\wedge^{k-(\ell-3)}(\mathcal{V}/\mathcal{W})).$$
There is a $\wedge^k(A)$-equivariant projection from the
image in $\wedge^k(\mathcal{V})$ of
$$\wedge^{\ell-3}({\mathcal W}) \otimes \wedge^{k-(\ell-3)}({\mathcal V})$$
onto $\wedge^{\ell-3}(\mathcal{W})\otimes 
\wedge^{k-(\ell-3)}(\mathcal{V}/\mathcal{W})$.
Therefore,
$0 \ne 
(\wedge^{k}(A)-1)^{\ell-1}(\wedge^{k}(\mathcal{V}))$.
\end{proof}

\begin{thm}
\label{unipex}
Suppose $\ell$ is a prime number, $\ell\ge 5$,
$V$ is a finite-di\-men\-sion\-al $\Q_\ell$-vector space, 
$T$ is a $\Z_\ell$-lattice in $V$,  
and $g$ is a quasi-unipotent linear operator on $V$ 
such that $g(T)=T$. 
Suppose $k$ and $m$ are positive integers 
and $2 \le k \le \dim(V)-2$.
If $(\wedge^k(g)-1)^{m(\ell-1)} \in \ell^m\End(\wedge^k(T))$, 
then $\wedge^k(g)$ is unipotent.  
\end{thm}

\begin{proof}
By Theorem \ref{localglobal},  
all the eigenvalues of 
$\wedge^k(g)$ are $\ell$-th roots of unity. 
In particular, $\wedge^k(g)^\ell$ is unipotent. 
By Lemma \ref{unipzl} 
there exists a $k$-th root of unity 
$\gamma \in \Z_\ell$ such that $(\gamma g)^\ell$ is unipotent. 
Replacing $g$ by $\gamma g$, we may assume
that $g^\ell$ is unipotent.

First, suppose that $g$ is semisimple. 
Then $g^\ell=1$ and $\wedge^k(g)^\ell=1$.
By Theorem 6.8 of \cite{serrelem} there is 
a $\wedge^k(g)$-invariant splitting of the free 
$\Z_\ell$-module $\wedge^k(T)$ into a direct sum of free 
$\Z_\ell$-modules
$\wedge^k(T)=P_1 \oplus P_2$
such that
$\wedge^k(g)$ acts as the identity on $P_1$, and
$$\wedge^k(g)^{\ell-1}+\cdots +\wedge^k(g)+1=0 \text{ on $P_2$.}$$
This implies easily that
$$(\wedge^k(g)-1)^{\ell-1}\in \ell \End(\wedge^k(T)).$$
If $g \ne 1$ then $g$ has an eigenvalue which is 
a primitive $\ell$-th root of unity, and 
therefore by Theorem \ref{localglobal},
$(g-1)^{\ell-2} \notin \ell \End(T)$.
If we let ${\mathcal V}=T/\ell T$, and let 
$A:{\mathcal V}\to{\mathcal V}$ 
be the linear operator induced by $g$, then
$(A-1)^{\ell-2} \ne 0$, but
$(\wedge^k(A)-1)^{\ell-1}=0$.
This contradicts Theorem \ref{unip}, and
proves that $g=1$ when $g$ is semisimple.

Next we will induct on the maximum of the multiplicities of the
roots of the minimal polynomial $P(t)$ of $g$ (i.e., on the
maximal size of the Jordan blocks for $g$).
Let $P_1(t)  \in \Z_\ell[t]$ be the monic polynomial whose roots
are the same as those of $P(t)$, but all with multiplicity one. 
Then $P_1$ divides $P$, and $P_1=P$ if and only 
if $g$ is semisimple. Let
$$T_0= \{x \in T \mid P_1(g)(x)=0\}.$$
Then $T_0$ is a pure free $\Z_\ell$-submodule of $T$ 
which is $g$-invariant, and the restriction $g_0:T_0\to T_0$ 
is semisimple. Let $T_1=T/T_0$ and let
$g_1$ denote the induced automorphism $g_1:T_1\to T_1$.
Then $T_1$ is a free $\Z_\ell$-module of finite rank, and 
the maximal multiplicity of a root of 
the minimal polynomial of
$$g_0\oplus g_1: T_0 \oplus T_1 \to T_0 \oplus T_1$$ 
is strictly less than that of $P(t)$, if $g$ is not
semisimple. 
Note that $(g_0\oplus g_1)^\ell$ is unipotent, 
since $g^\ell$ is unipotent.
Further, $g_0\oplus g_1$ is unipotent 
if and only if $g$ is unipotent.
To apply induction and finish the proof, 
it suffices to check that
\begin{equation}
\label{tteqn}
(\wedge^k(g_0\oplus g_1)-1)^{m(\ell-1)} 
\in \ell^m \End(\wedge^k(T_0\oplus T_1)).
\end{equation}
For $0\le i \le k$,
let $H_i$ be the image of the natural homomorphism 
$$\wedge^i(T_0)\otimes \wedge^{k-i}(T_1) \to \wedge^k(T).$$
Then 
$$H_i/H_{i+1} \cong \wedge^i(T_0)\otimes \wedge^{k-i}(T_1),$$
and 
$$\wedge^k(T)=H_0 \supset H_1 \supset \cdots \supset 
H_k=\wedge^k(T_0)\supset H_{k+1}:=0$$
is a natural filtration of $\wedge^k(g)$-stable pure 
$\Z_\ell$-submodules of $\wedge^k(T)$.
Since $H_i$ is pure in $\wedge^k(T)$, we have
$H_i\cap\ell^m\wedge^k(T)=\ell^mH_i$. Therefore,
$$(\wedge^k(g)-1)^{m(\ell-1)}(H_i) \subseteq \ell^mH_i.$$ 
Since
$$\oplus_{i=0}^k(H_i/H_{i+1}) =
\oplus_{i=0}^k(\wedge^i(T_0)\otimes \wedge^{k-i}(T_1))=
\wedge^k(T_0\oplus T_1),$$
we have (\ref{tteqn}).
\end{proof}

\section{Abelian variety lemmas}
\label{lemsect}

As stated earlier,
we suppose $X$ is an abelian variety over a field $F$, 
$v$ is a discrete valuation on $F$ of residue characteristic 
$p \ge 0$, and $\ell$ is a prime different from $p$.
Recall that $\I$ is the inertia subgroup at 
${\bar v}$ of $\Gal(F^s/F)$.

\begin{thm}[Galois Criterion for Semistability]
\label{galcrit}
The following are equivalent:
\begin{enumerate}
\item[(i)] $X$ is semistable at $v$,
\item[(ii)] $\I$ acts unipotently on $T_\ell(X)$; i.e.,
all the eigenvalues of $\rho_{\ell,X}(\sigma)$ are $1$  
for every $\sigma \in \I$,
\item[(iii)] for every $\sigma \in \I$, 
$(\rho_{\ell,X}(\sigma) - 1)^2 = 0$.
\end{enumerate}
\end{thm}

\begin{proof}
See 3.5 and 3.8 of \cite{SGA7}, 
and Theorem 6 on p.~184 of \cite{BLR}.
\end{proof}

\begin{lem}
\label{level2}
Suppose $\sigma \in \I$.
Then: 
\begin{enumerate}
\item[{(i)}] $\rho_{\ell,X}(\sigma)$ is unipotent 
if and only if  $(\rho_{\ell,X}(\sigma)-1)^2=0$;
\item[{(ii)}] $-\rho_{\ell,X}(\sigma)$ is unipotent 
if and only if  $(\rho_{\ell,X}(\sigma)+1)^2=0$.
\end{enumerate}
\end{lem}

\begin{proof}
There exists a finite Galois extension 
$L \subset F^s$ of $F$ such
that if $w$ is the restriction of $\bar{v}$ to $L$ then 
$X$ is semistable at $w$
(see Prop.~3.6 of \cite{SGA7}).
Let $\I_w=\I\cap\Gal(F^s/L)$ be the corresponding inertia group,
let $m=[L:F]$, and let $g=\rho_{\ell,X}(\sigma)$. 
Then $\sigma^m \in \I_w$. 
By Theorem \ref{galcrit}, $(g^m-1)^2=0$. 
Now apply Lemma \ref{level2lem}.
\end{proof}

The following result follows immediately from 
Lemmas \ref{lin}iii and \ref{level2}.

\begin{prop}
\label{tatemodlem}
Suppose $\sigma \in \I$.
The following are equivalent:
\begin{enumerate}
\item[{(i)}] there exists $\gamma_\sigma \in \Q_\ell$ such that
$\gamma_\sigma\rho_{\ell,X}(\sigma)$ is unipotent,
\item[{(ii)}] either $\rho_{\ell,X}(\sigma)$ or 
$-\rho_{\ell,X}(\sigma)$ is unipotent,
\item[{(iii)}] 
either $(\rho_{\ell,X}(\sigma)-1)^2=0$ or $(\rho_{\ell,X}(\sigma)+1)^2=0$.
\end{enumerate}
\end{prop}

\begin{prop}
\label{padd}
Let $\J \subset \I$ denote the first ramification group,
and let $\tau$ be a lift to $\I$ of a topological generator
of the procyclic group $\I/\J$. 
The following are equivalent:
\begin{enumerate}
\item[(i)] 
$X$ has purely additive reduction at $v$,
\item[(ii)] 
$1$ is not an eigenvalue for the action of $\tau$ on $V_{\ell}(X)^\J$,
\item[(iii)] 
$V_{\ell}(X)^{\I}=0$.
\end{enumerate}
\end{prop}

\begin{proof}
The equivalence of (ii) and (iii) is obvious. 
For the equivalence of (i) and (ii) see 
Corollary 1.10 of \cite{LenstraOort}.
\end{proof}

\section{Higher cohomology groups of abelian varieties}
\label{highcohsect}

Write $V = V_\ell(X)$ and $T = T_\ell(X)$, and recall that
$H^k_\ell = H^k_{\text{\'et}}({\bar X},\Z_\ell)$.
The image of $\rho_{\ell,X}$ lies in the
symplectic group $\Sp(V)$, by the Galois-equivariance of the 
Weil pairing, and the fact that the inertia group acts as the 
identity on the $\ell$-power roots of unity. 

\begin{asmpt}
\label{assumption}
For the remainder of this paper (except for Remark \ref{henrem})
we will assume that if $p=2$ then the valuation ring
is henselian.
\end{asmpt}

\begin{defn}
\label{bu1}
If $p \ne 2$ then
we say that $X$ is {\em briefly unstable} at $v$ if 
$X$ is purely additive at $v$ and 
becomes semistable above $v$ over a quadratic 
separable extension of $F$. 
\end{defn}

\begin{defn}
\label{bu2}
If $p=2$ then we say that $X$ is {\em briefly unstable} at $v$ if
$X$ is purely additive at $v$ and there exists a finite 
unramified extension $M$ of $F$ such that $X$ is semistable 
above $v$ over a quadratic separable extension of $M$. 
\end{defn}

\begin{rem}
By Theorem \ref{galcrit}, the quadratic extension in
Definitions \ref{bu1} and \ref{bu2} is ramified over $v$.
\end{rem}

\begin{thm}
\label{tatemod}
Suppose $X$ is an abelian variety over a field $F$, 
$v$ is a discrete valuation on $F$ of residue characteristic 
$p \ge 0$, and $\ell$ is a prime different from $p$. 
Then the following are equivalent:
\begin{enumerate}
\item[{(a)}] $X$ is either semistable or briefly unstable at $v$, 
\item[{(b)}] for each $\sigma \in \I$,
either $\rho_{\ell,X}(\sigma)$ or $-\rho_{\ell,X}(\sigma)$ is unipotent.
\end{enumerate}
\end{thm}

\begin{proof}
Assume (a) holds.
By Theorem \ref{galcrit}, we may reduce to the case where 
$X$ has purely additive reduction at $v$, $M$ is a finite 
unramified extension of $F$, and $L$ is a 
quadratic separable extension of $M$ over which $X$ is 
semistable above $v$. 
Then by Theorem \ref{galcrit}, 
$\rho_{\ell,X}(\sigma)^2$ is unipotent for all $\sigma \in \I$.
Let $\J \subset \I$ be the first ramification subgroup. 
Then $\J$, and therefore $\rho_{\ell,X}(\J)$, 
is either trivial (if $p=0$) or a pro-$p$-group. 
By \cite{LenstraOort} (see pp.~282--283), 
since $\ell \ne p$, $\rho_{\ell,X}(\J)$ is either trivial or a finite
$p$-group.  
If $s \in \rho_{\ell,X}(\J) \subset
\rho_{\ell,X}(\I)$, then $s^2$ is unipotent and 
has finite order, and thus $s^2=1$.
It follows that either
$\rho_{\ell,X}(\J)=\{1\}$, 
or $p=2$ and $\rho_{\ell,X}(\J)$ is a finite commutative 
group of exponent $2$.
 
Suppose that $\rho_{\ell,X}(\J)=\{1\}$. 
Then $V^{\J}=V$.
Let $\tau$ be a lift to $\I$ of a topological generator
of the procyclic group $\I/\J$. 
Then $g:=\rho_{\ell,X}(\tau)$ generates the
procyclic group $\rho_{\ell,X}(\I)$.
By Prop.\ \ref{padd}, $1$ is not an eigenvalue of $g$. 
Since $g^2$ is unipotent, the only eigenvalue of $g$ is $-1$, 
i.e., $-g$ is unipotent. For each integer $i$ either
$g^i$ or $-g^i$ is unipotent. 
Since in the $\ell$-adic topology
the set of integral powers of $g$ is
dense in $\rho_{\ell,X}(\I)$ and
the set of unipotent operators in
$\Aut(T)$ is closed, therefore
for each $\sigma \in \I$ either $\rho_{\ell,X}(\sigma)$ or
$-\rho_{\ell,X}(\sigma)$ is unipotent.
 
We may thus assume that $\rho_{\ell,X}(\J)\ne \{1\}$, $p=2$, and
$\rho_{\ell,X}(\J)$ is a finite commutative group of exponent $2$. 
We may assume that $L \subset F^{s}$.
Let $w$ be the restriction of $\bar{v}$ to $L$. 
Let $\I_{w}$ denote the inertia subgroup at ${\bar v}$ of $\Gal(F^{s}/L)$. 
Clearly,  
$\J_w:=\J \cap \I_w$ is 
the first ramification subgroup of $\I_w$,
and $\J_w$ has index $2$ in $\J$.
Since $X$ is semistable at $w$, $\rho_{\ell,X}(\sigma)$ is
unipotent for all $\sigma \in \I_w$.
Since $\rho_{\ell,X}(\J)$ is finite,
$\rho_{\ell,X}(\sigma)=1$ for all $\sigma \in \J_w$.
Since $\J_w$ has index $2$ in $\J$, therefore
$\rho_{\ell,X}(\J)$ has order $2$.
 
Since $L/F$ is wild quadratic, therefore 
the inclusion $\I_w \subset \I$ induces a natural isomorphism
$\I_w/\J_w = \I/\J$. 
Let $\tau_w$ be a lift to $\I_w$ 
of a topological generator of $\I_w/\J_w=\I/\J$. 
Since $\rho_{\ell,X}(\tau_w)$ is unipotent, and 
$X$ has purely additive reduction at $v$, therefore $V^{\J}=0$
by Prop.\ \ref{padd}.

Clearly $\J_w$ is normal in $\I$,
since it is the intersection of normal subgroups.  
We can view $\rho_{\ell,X}$ as a homomorphism from $\I/\J_w$
to $\Aut(T)$.
The image of 
$$\I_w \to \I_w/\J_w \subset \I/\J_w =\I/\J \times \I/\I_w$$
is $\I/\J \times \{1\}$. Therefore
$\rho_{\ell,X}(\I/\J \times \{1\})$ consists of unipotent operators.
 Let $s$ be the non-trivial element of $\I/\I_w$ and 
let $h=\rho_{\ell,X}(1 \times s)$. Then $h^2=1$. If $h=1$ then
$\rho_{\ell,X}(\I/\J_w)$ consists of unipotent operators, so 
$X$ is semistable at $v$, which is not the case. 
So $h \ne 1$. 
If $h=-1$ then 
$\rho_{\ell,X}(\I/\J_w)$ is the union of 
$\rho_{\ell,X}(\I/\J \times \{1\})$ and 
$-\rho_{\ell,X}(\I/\J \times \{1\})$. 
Therefore for each $g \in \rho_{\ell,X}(\I/\J_w)$, either $g$ or $-g$ is
unipotent. So we have reduced to the case where $h \ne \pm 1$. 
But then $V^\J$,
the eigenspace of $h$ corresponding to the eigenvalue $1$,
is non-zero.
This contradiction proves that (a) implies (b).

To prove that (b) implies (a), 
suppose that $X$ is not semistable at $v$,
and suppose that for each $\sigma \in \I$
either $\rho_{\ell,X}(\sigma)$ or $-\rho_{\ell,X}(\sigma)$ is unipotent.
By Theorem \ref{galcrit}, $\rho_{\ell,X}(\sigma)$ is not unipotent for
some $\sigma\in\I$. For such a $\sigma$, the eigenvalues of 
$\rho_{\ell,X}(\sigma)$ are all $-1$. Thus $V^{\I}=0$.
By Prop.\ \ref{padd}, $X$ has purely additive reduction at $v$. 
Let 
\begin{equation}
\label{Ivxeqn}
\I_{v,X}=\{\sigma \in \I\mid \rho_{\ell,X}(\sigma) \text{ is unipotent} \}.
\end{equation}
Then $\I_{v,X} \ne \I$.
It is known  (see pp.\ 354--355 of \cite{SGA7}
and \S 4 of \cite{inertia}) that 
$\I_{v,X}$ is an open normal subgroup of finite index in $\I$.
Since $\rho_{\ell,X}(I_{v,X})$ consists of unipotent operators, 
$V^{I_{v,X}} \ne 0$
by a theorem of Kolchin (p.~35 of \cite{SerreLie}).
The restriction map
$\rho':I \to \Aut(V^{\I_{v,X}})$ factors through the
finite group $I/I_{v,X}$. 
Therefore the image of $\rho'$ is finite.
If $\sigma \in \I-\I_{v,X}$, then $-\rho'(\sigma)$ is
unipotent and of finite order, and thus $\rho'(\sigma)=-1$
on $V^{\I_{v,X}}$. 
Therefore $\rho'$ has kernel $\I_{v,X}$ and image $\{\pm 1\}$,
so $[\I:\I_{v,X}]=2$.
   
First assume $p \ne 2$. 
Then $\I$ has exactly one subgroup of index $2$ and thus
this subgroup must be $\I_{v,X}$.
Let $L/F$ be a ramified separable quadratic extension. 
We may assume that $L\subset F^s$. 
Let $w$ be the restriction of $\bar{v}$ to $L$. The
corresponding inertia group $\I_w \subset \Gal(F^s/L)$ 
has index $2$ in $\I$ and therefore is $\I_{v,X}$.
By Theorem \ref{galcrit}, $X$ is semistable at $w$. 
  
Now assume that $p=2$. 
Let $F^{ur}=(F^s)^\I$, the maximal extension of $F$ 
unramified above $v$. 
The valuation ring of $F^{ur}$ is henselian
and the residue field is separably closed. 
Let $L$ be the quadratic extension of $F^{ur}$ 
corresponding to $\I_{v,X} \subset \I\cong\Gal(F^s/F^{ur})$. 
Then $L=F^{ur}(\sqrt{c})$ for some $c\in F^{ur}$,
$F(c)$ is unramified above $v$ over $F$,
and $X$ is semistable over the ramified quadratic extension
$F(\sqrt{c})$ of $F(c)$.
\end{proof}

\begin{thm}
\label{cohcor}
Suppose $X$ is an abelian variety over a field $F$, 
suppose $v$ is a discrete valuation on $F$ of residue 
characteristic $p \ge 0$, suppose $k$ is a positive integer, 
suppose $k < 2\dim(X)$,
and suppose $\ell$ is a prime number
not equal to $p$.
\begin{enumerate}
\item[{(i)}] If $k$ is odd then the following are
equivalent:
\begin{enumerate}
\item[{(a)}]  $X$ is semistable at $v$,
\item[{(b)}]  for each $\sigma\in\I$,
$\rho_{\ell,k}(\sigma)$ is unipotent, 
\item[{(c)}]  for each $\sigma\in\I$,
$(\rho_{\ell,k}(\sigma)-1)^{k+1}=0$. 
\end{enumerate}
\item[{(ii)}] If $k$ is even 
then the following are equivalent:
\begin{enumerate}
\item[{(a)}]  $X$ is either semistable or briefly unstable at $v$,
\item[{(b)}]  for each $\sigma\in\I$,
$\rho_{\ell,k}(\sigma)$ is unipotent, 
\item[{(c)}]  for each $\sigma\in\I$,
$(\rho_{\ell,k}(\sigma)-1)^{k+1}=0$.
\end{enumerate}
\item[{(iii)}] If $k$ is odd  
then the following are equivalent:
\begin{enumerate}
\item[{(a)}]  $X$ is either semistable or briefly unstable at $v$,
\item[{(b)}]  for each $\sigma\in\I$, 
either $\rho_{\ell,k}(\sigma)$ or $-\rho_{\ell,k}(\sigma)$ 
is unipotent, 
\item[{(c)}]  for each $\sigma\in\I$, 
either $(\rho_{\ell,k}(\sigma)-1)^{k+1}=0$ or
 $(\rho_{\ell,k}(\sigma)+1)^{k+1}=0$. 
\end{enumerate}
\end{enumerate}
\end{thm}

\begin{proof}
Clearly, (c) implies (b). That (a) implies (c) follows from 
Remark \ref{astrem}, combined with 
Theorem \ref{galcrit} and Lemmas \ref{level2}i and \ref{lin}i
for (i),  with 
Theorem \ref{tatemod} and Lemma \ref{lin}ii for (ii),
and with
Theorem \ref{tatemod} and Lemmas \ref{level2} and \ref{lin}i for (iii).
Suppose we have (b). 
To conclude (a),
apply Lemma \ref{lin}iv and Remark \ref{astrem},
combined with 
Theorem \ref{tatemod} for (ii) and (iii) 
and with Theorem \ref{galcrit} for (i).
\end{proof}

\begin{cor}
\label{highercoha}
Suppose $X$ is an abelian variety over a field $F$, 
and $v$ is a discrete 
valuation on $F$ of residue characteristic $p \ge 0$.
Suppose $k$ and $r$ are positive integers, 
$k < 2\dim(X)$, and $k<r$. Suppose $\sigma \in \I$. 
\begin{enumerate}
\item[{(i)}] 
If either $X$ is semistable at $v$, 
or $k$ is even and $X$ is briefly unstable at $v$,  
then 
$$(\sigma - 1)^{r}H^k_{\text{\'et}}({\bar X}, \Z_{\ell}) = 0$$ 
for every prime $\ell \ne p$,
and
$$(\sigma - 1)^{r}H^k_{\text{\'et}}({\bar X}, \Z/n\Z) = 0$$ 
for every positive integer $n$ not divisible by $p$.
\item[{(ii)}] 
If $k$ is odd and $X$ is briefly unstable at $v$, then 
for every prime $\ell \ne p$,
either
$$(\sigma - 1)^{r}H^k_{\text{\'et}}({\bar X}, \Z_{\ell}) = 0
\quad \text{ or }\quad 
(\sigma + 1)^{r}H^k_{\text{\'et}}({\bar X}, \Z_{\ell}) = 0,$$
and for every positive integer $n$ not divisible by $p$,
either
$$(\sigma - 1)^{r}H^k_{\text{\'et}}({\bar X}, \Z/n\Z) = 0
\quad \text{ or }\quad 
(\sigma + 1)^{r}H^k_{\text{\'et}}({\bar X}, \Z/n\Z) = 0.$$  
\end{enumerate}
\end{cor}

\begin{proof}
The first parts follow from Theorem \ref{cohcor}, 
since $r\ge k+1$.
The second parts follow from the first parts for all
prime divisors $\ell$ of $n$, since for all $i$,
$$H^k_{\text{\'et}}({\bar X}, \Z/\ell^i\Z)=
H^k_{\text{\'et}}({\bar X}, \Z_\ell)\otimes \Z/\ell^i\Z.$$
\end{proof}

\begin{thm}
\label{highercoh}
Suppose $X$ is an abelian variety over a field $F$, 
and $v$ is a discrete 
valuation on $F$ of residue characteristic $p \ge 0$.
Suppose $k$, $n$, and $r$ are positive integers, 
$k < 2\dim(X)$, $n$ is not divisible by $p$, 
and $n \notin N(r)$.
\begin{enumerate}
\item[{(i)}] 
Suppose that  
$(\sigma - 1)^{r}H^k_{\text{\'et}}({\bar X}, \Z/n\Z) = 0$ 
for all $\sigma \in \I$.
Then either $X$ is semistable at $v$,
or $k$ is even and $X$ is briefly unstable at $v$.
\item[{(ii)}] 
Suppose $k$ is odd, and suppose that for each $\sigma \in \I$ either 
$$(\sigma - 1)^{r}H^k_{\text{\'et}}({\bar X}, \Z/n\Z) = 0 
\quad \text{ or }\quad
(\sigma + 1)^{r}H^k_{\text{\'et}}({\bar X}, \Z/n\Z) = 0.$$
Then either $X$ is semistable at $v$,
or $X$ is briefly unstable at $v$.
\end{enumerate}
\end{thm}

\begin{proof}
Recall (\cite{SGA7}, Thm.~4.3) that 
the characteristic polynomial of $\rho_{\ell,X}(\sigma)$ has 
integer coefficients and does not depend on the choice of 
$\ell \ne p$. 
Therefore the characteristic polynomial 
$\P_{\sigma}$ of $(\rho_{\ell,k}(\sigma)-1)^r/n$ 
has coefficients in $\Z[1/n]$ and does not
depend on the choice of $\ell \ne p$.
Suppose that 
$(\sigma - 1)^{r}H^k_{\text{\'et}}({\bar X}, \Z/n\Z) = 0$.  
Then for all prime divisors $\ell$ of $n$, 
$(\rho_{\ell,k}(\sigma)-1)^rH_\ell^k \subseteq nH_\ell^k$, 
so $\P_{\sigma}$ has coefficients in $\Z_{\ell}$.
Thus $\P_{\sigma}$ has integer coefficients.
Since $X$ is semistable over a finite separable 
extension of $F$,
by Theorem \ref{galcrit}, Lemma \ref{lin}i,ii, and Theorem \ref{cohcor}
there is a positive integer $m$ such that
$(\rho_{\ell,k}(\sigma^m) - 1)^{k+1} = 0$.
Let $\alpha$ be an eigenvalue of $\rho_{\ell,k}(\sigma)$. 
Then 
$(\alpha^m - 1)^{k+1} = 0$, so $\alpha^m = 1$. 
Since $(\alpha-1)^r/n$ is an eigenvalue of
$(\rho_{\ell,k}(\sigma)-1)^r/n$ and therefore is a root of $\P_{\sigma}$,
it is an algebraic integer.
Since $n \notin N(r)$,
we have $\alpha = 1$ by Theorem \ref{localglobal2}.
Thus $\rho_{\ell,k}(\sigma)$ is unipotent. 
Applying Theorem \ref{cohcor}, we have (i).
To obtain (ii), replace $-1$ by $+1$ in the above argument.
\end{proof}

Next we specialize Theorem \ref{highercoh}i to the case $r=1$.
Note that when $k=1$ we recover Raynaud's criterion for
semistability (Prop~4.7 of \cite{SGA7}).

\begin{thm}
\label{raynaudgen}
Suppose $X$ is an abelian variety over a field $F$, 
$v$ is a discrete valuation on $F$, 
$k$ and $n$ are integers, 
$0 < k < 2\dim(X)$, $n \ge 3$, 
$n$ is not divisible by the residue characteristic, and
$\I$ acts as the identity on 
$H^k_{\text{\'et}}({\bar X},\Z/n\Z)$.
Then either $X$ is semistable at $v$, or
$k$ is even  
and $X$ is briefly unstable at $v$.
\end{thm}

The next result is an immediate corollary of 
Theorems \ref{highercoh} and  \ref{cohcor}
and Cor.\ \ref{highercoha}.

\begin{thm}
\label{highercohcor}
Suppose $X$ is an abelian variety over a field $F$, 
suppose $v$ is a discrete valuation on $F$ of residue 
characteristic $p \ge 0$, suppose $k$, $n$, and $r$
are positive integers, and suppose $\ell$ is a prime number.
Suppose $k < 2\dim(X)$, 
suppose $k<r$, suppose $n\ell$ is not divisible by $p$, 
and suppose $n \notin N(r)$.
\begin{enumerate}
\item[{(i)}] If $k$ is odd, then the following are
equivalent:
\begin{enumerate}
\item[{(a)}]  $X$ is semistable at $v$,
\item[{(b)}] for each $\sigma \in \I$,
$(\sigma - 1)^{r}H^k_{\text{\'et}}({\bar X}, 
\Z_{\ell}) = 0$, 
\item[{(c)}]  for each $\sigma \in \I$,
$(\sigma - 1)^{r}H^k_{\text{\'et}}({\bar X},\Z/n\Z) = 0$. 
\end{enumerate}
\item[{(ii)}] If $k$ is even 
then the following are equivalent:
\begin{enumerate}
\item[{(a)}] $X$ is either semistable or briefly unstable at $v$,
\item[{(b)}] for each $\sigma \in \I$,
$(\sigma - 1)^{r}H^k_{\text{\'et}}({\bar X},\Z_{\ell}) = 
0$,
\item[{(c)}] for each $\sigma \in \I$,
$(\sigma - 1)^{r}H^k_{\text{\'et}}({\bar X},\Z/n\Z) = 0$. 
\end{enumerate}
\item[{(iii)}] If $k$ is odd 
then the following are equivalent:
\begin{enumerate}
\item[{(a)}] $X$ is either semistable or briefly unstable at $v$,
\item[{(b)}] for each $\sigma \in \I$ either 
$(\sigma - 1)^{r}$ 
or
$(\sigma + 1)^{r}$ 
kills $H^k_{\text{\'et}}({\bar X},\Z_{\ell})$. 
\item[{(c)}] for each $\sigma \in \I$ either
$(\sigma - 1)^{r}$ 
or
$(\sigma + 1)^{r}$ 
kills $H^k_{\text{\'et}}({\bar X}, \Z/n\Z)$. 
\end{enumerate}
\end{enumerate}
\end{thm}

\section{Non-extremal cohomology groups}
\label{midk}

The results of \S\ref{highcohsect} are sharp when $k = 1$ or
$2d-1$.
In this section we obtain sharp results when $1 < k < 2d-1$.
Let 
$$N'(r)=
 \{\ell^m \mid 0 \le m(\ell - 1) < r, \text{
or $\ell=2$ or $3$ and $r=m(\ell-1)$\}.}$$
Clearly, $N'(r) \subseteq N(r)$, and a prime power $\ell^m$
lies in $N(r) - N'(r)$ 
if and only if $\ell \ge 5$ and $r=m(\ell-1)$.

\begin{thm}
\label{extreme}
Suppose $X$ is an abelian variety over a field $F$, 
suppose $v$ is a discrete valuation on $F$ of residue 
characteristic $p \ge 0$, suppose $k$ and $r$
are positive integers, and suppose $n$ is a positive integer
which is not divisible by $p$.
Suppose $2 \le k \le 2\dim(X)-2$, suppose $k<r$, and suppose
$n \notin N'(r)$.
\begin{enumerate}
\item[{(i)}] If $k$ is odd, then the following are
equivalent:
\begin{enumerate}
\item[{(a)}]  $X$ is semistable at $v$,
\item[{(b)}] for every $\sigma \in \I$,
$(\sigma - 1)^{r}H^k_{\text{\'et}}({\bar X}, \Z/n\Z) = 0$. 
\end{enumerate}
\item[{(ii)}] If $k$ is even 
then the following are equivalent:
\begin{enumerate}
\item[{(a)}]  $X$ is either semistable or briefly unstable at $v$,
\item[{(b)}] for every $\sigma \in \I$,
$(\sigma - 1)^{r}H^k_{\text{\'et}}({\bar X}, \Z/n\Z) = 0$. 
\end{enumerate}
\item[{(iii)}] If $k$ is odd 
then the following are equivalent:
\begin{enumerate}
\item[{(a)}]  $X$ is either semistable or briefly unstable at $v$,
\item[{(b)}] for each $\sigma \in \I$ either 
$(\sigma - 1)^{r}$ 
or
$(\sigma + 1)^{r}$
kills $H^k_{\text{\'et}}({\bar X}, \Z/n\Z)$. 
\end{enumerate}
\end{enumerate}
\end{thm}

\begin{proof}
If  $n \notin N(r)$ then the assertion is contained in
Theorem \ref{highercohcor}.
Thus we may assume that $n=\ell^m$ with $\ell \ge 5$ and $r=m(\ell-1)$.
If (a) holds, then (b) follows from Corollary \ref{highercoha}.
If (b) holds, then (a) follows from 
Theorems \ref{unipex} and \ref{cohcor}.
\end{proof}

The following example shows that the condition 
$n \notin N'(r)$ is sharp. 

\begin{ex}
\label{example}
Suppose that $F$ is a discrete valuation field,
$m$ and $r$ are positive integers, $\ell$
is a prime number, and $n=\ell^m$. 
Suppose that either $m(\ell-1)<r$, or
$\ell=2$ and $r=m=m(\ell-1)$, or
$\ell=3$ and $r=2m=m(\ell-1)$.
Let $X_1$, $Y$ be abelian varieties of positive dimension
with good reduction over $F$, and such that $Y$ has an
automorphism of exact order $\ell$.
Let $L$ be a totally ramified
degree $\ell$ extension of $F$, and let $X_2$ be the $L/F$-form
of $Y$ corresponding to a character
$\Gal(L/F) \hookrightarrow \Aut(Y)$.
Then $X=X_1 \times X_2$ is neither semistable nor 
purely additive.
However, by Lemma \ref{groupring} we have 
$(\rho_{\ell,k}(\sigma)-1)^{r} \in n \End(H_\ell^k)$
for every $\sigma\in\I$.
\end{ex}

\section{Semistability over quadratic extensions}
\label{rmksect}

\begin{rem}
\label{henrem}
When the valuation ring is henselian, then
$v$ extends uniquely to $F^s$, $\I$ is normal in $\Gal(F^s/F)$, 
and the field $(F^s)^\I$ is the (unique) 
maximal extension of $F$ unramified over $v$.
Note that 
Assumption \ref{assumption} can be dropped in 
Theorems \ref{cohcor}i, \ref{highercohcor}i, and \ref{extreme}i,
in Corollary \ref{highercoha}i if $X$ is semistable,
and in Theorems \ref{highercoh}i and \ref{raynaudgen} 
when $k$ is odd, since it is not used.
(Assumption \ref{assumption} is only used
when we deal with brief instability.)
\end{rem}

\begin{lem}
\label{2lem}
Suppose $X$ is an abelian variety over a field $F$, 
$v$ is a discrete valuation on $F$ 
of residue characteristic $p\ge 0$, and
the valuation ring is henselian.
Suppose there exists a finite unramified extension $M$ of $F$ 
such that $X$ is semistable above $v$ over a quadratic separable 
extension $L$ of $M$.
Suppose that if $p=2$, then the residue field either is 
separably closed or is algebraic over $\F_2$ (e.g., is finite).
Then $X$ is semistable above $v$ over a quadratic 
separable extension of $F$.
\end{lem}

\begin{proof}
If $K$ is a subfield of $F^s$, let $\I_K$ denote 
the inertia subgroup of $\I$ corresponding to $K$.
Note that $\I_L$ is an open subgroup of $\I$ of index $2$.

When $p\ne 2$, then $F$ has a tamely ramified 
separable quadratic extension $L'$.   
Since $p \ne 2$, $\I$ has 
exactly one open subgroup of index $2$, 
so $\I_L=\I_{L'}$. 
By Theorem \ref{galcrit} over $L$ and over $L'$, 
$X$ is semistable above $v$ over $L'$.
Now suppose $p=2$. 

If the residue field is separably closed, then 
$F$ has no non-trivial unramified extensions, so $M=F$. 

Suppose the residue field $k$ is algebraic over $\F_2$, 
and let $G_k:=\Gal(k^s/k)$. 
Then $G_k$ is a torsion-free procyclic group,
since it is a closed subgroup of $\hat{\Z}$.
Let $F^{ur}$ be the maximal unramified extension of $F$. Then
$\I=\Gal(F^s/F^{ur})$.
Since the valuation ring is henselian, $G_k=\Gal(F^{ur}/F)$.
We may assume that $X$ is not semistable at $v$.
Let $\I_{v,X}$ be the group defined by formula (\ref{Ivxeqn})
(proof of Theorem \ref{tatemod}) with $\ell=3$.
Applying Theorem \ref{galcrit} over $F$ and over $L$ (for 
$\ell=3$) shows that $\I_{v,X}$ is a proper subgroup of $\I$ and
$\I_L \subseteq \I_{v,X} \subset \I$.
Since $[\I:\I_L]=2$, we have $\I_{v,X}=I_L$ and $[\I:\I_{v,X}]=2$.
Let $F'$ be the quadratic extension of $F^{ur}$ 
cut out by $\I_{v,X}$. 
Then $F'/F$ is Galois, since the group $\I_{v,X}$ is the 
intersection of $\I$ and 
$$\{\sigma \in \Gal(F^s/F)\mid \rho_{3,X}(\sigma) 
\text{ is unipotent} \},$$
and both are stable under conjugation by $\Gal(F^s/F)$.
Further, $X$ is semistable (above $v$) over $F'$ by
Theorem \ref{galcrit}, and
$G:=\Gal(F'/F)$ is an extension of $G_k$ by $C:=\Gal(F'/F^{ur})$.
This extension is central, 
since $C$ has order $2$ and thus has no non-trivial automorphisms. 
Since every group whose quotient by its center is 
(pro)cyclic must be commutative, $G$ is commutative.
Let $\Delta$ be the subset of squares in $G$.
Then $\Delta$ is a closed subgroup, $\Delta\cap C=1$, and
$[G:\Delta]=2$ or $4$ 
($G_k$ is procyclic, so $2G_k$ has index $1$ or $2$ in $G_k$). 
Thus $\Delta$ is open and $G/\Delta$ either 
has order $2$ (in which case let $H=\Delta$)
or is a product of two groups of order $2$,
one of which is the image of $C$
(in which case let $H$ be the preimage in $G$ 
of the other one).
Then $H$ is an open subgroup of $G$ of index $2$,  
so the corresponding subfield $L':=(F')^H$ is quadratic over $F$.
Since $H\cap C=1$, therefore $F'/L'$ is unramified and 
so $X$ is semistable over $L'$. 
\end{proof}

\begin{cor}
Suppose that if the residue characteristic is $2$, then
the residue field either is 
separably closed or is algebraic over $\F_2$. Then
we may replace ``briefly unstable at $v$'' by 
``purely additive at $v$ and becomes semistable 
above $v$ over a quadratic 
separable extension of $F$'' in the results of 
\S\S\ref{highcohsect}--\ref{midk}.
\end{cor}

\begin{rem}
The proof of Lemma \ref{2lem} shows that the condition
that the residue field $k$ be algebraic over $\F_2$
can be replaced by the condition that
$\Gal(k^s/k)$ be a torsion-free procyclic group.
\end{rem}

\begin{rem}
The group $\I_{v,X}$ defined in (\ref{Ivxeqn}) is 
independent of $\ell$ (see p.~355
of \cite{SGA7} and
Theorem 4.2 of \cite{inertia}).
It follows that 
for each fixed $\sigma\in\I$,
whether or not $(\sigma-1)^r$ kills 
$H^k_{\text{\'et}}({\bar X},\Z_{\ell})$
(or $H^k_{\text{\'et}}({\bar X},\Z/n\Z)$ for $n$ not a power of $2$)
is independent of $\ell$ (and $n$), and 
depends only on whether or not $\sigma\in\I_{v,X}$.
\end{rem}

\end{document}